\documentclass[12pt]{article}
\usepackage{graphicx}
\usepackage{amssymb}
\usepackage{amsfonts}
\usepackage{amsthm}
\usepackage{amsfonts}

\newcommand \ra {\rightarrow}

\newcommand{\ba}[1]{\begin{array}{#1}}
\newcommand{\ea}{\end{array}}

\newcommand{\be}{\begin{equation}}
\newcommand{\ee}{\end{equation}}
\newcommand{\bea}{\begin{eqnarray}}
\newcommand{\eea}{\end{eqnarray}}
\newcommand{\beann}{\begin{eqnarray*}}
\newcommand{\eeann}{\end{eqnarray*}}

\def\reff#1{(\ref{#1})}
\newtheorem{conjecture}{Conjecture}

\bibstyle{ams}

\begin{document}

\title{Conformal invariance of \\ the loop-erased percolation explorer}

\author{Tom Kennedy
\\Department of Mathematics
\\University of Arizona
\\Tucson, AZ 85721
\\ email: tgk@math.arizona.edu
}

\maketitle 

\begin{abstract}

We consider critical percolation on the triangular lattice in a bounded 
simply connected domain with boundary conditions that force an interface 
between two prescribed boundary points.
We say the interface forms a ``near-loop'' when it comes within one lattice
spacing of itself. We define a new curve by erasing these near-loops as 
we traverse the interface. Our Monte Carlo simulations of this model lead
us to conclude that the scaling limit of this loop-erased percolation interface
is conformally invariant and has fractal dimension $4/3$. However, it is not 
SLE$_{8/3}$. We also consider the process in which a near-loop is 
when the explorer comes within two lattice spacings of itself. 
\end{abstract}

\section{Introduction}

The scaling limits of many two-dimensional models from statistical mechanics
are conformally invariant when the model is critical. Many of 
these conformally invariant scaling limits are described by the 
Schramm-Loewner evolution (SLE$_\kappa$) for some value of the parameter 
$\kappa$ \cite{schramm_sle}. (For an exposition of SLE see \cite{lawler_book}.)
In this paper we introduce a new stochastic process arising in critical 
percolation, an important example of such models.
(For an exposition of percolation and its relation to SLE see
\cite {werner2007lectures}.)
We focus on critical site percolation on the triangular
lattice, the one case in which the conformal invariance has been proven. 

We can think of the sites in the triangular lattice
as the centers of the hexagons in a hexagonal lattice, and we will define
our model using the hexagonal lattice. 
Let $D$ be a bounded simply connected domain. Fix two sites $z$ and $w$ 
on its boundary. For a lattice spacing $\delta$ we let $D_\delta$ be 
a collection of hexagons which approximates $D$ and let $z_\delta$ and 
$w_\delta$ be sites in the hexagonal lattice which approximate $z$ and $w$. 
We color the hexagons along the boundary of $D_\delta$ going from $z_\delta$ to 
$w_\delta$ in the clockwise direction white,
and then color the hexagons along the 
boundary from $w_\delta$ back to $z_\delta$ in the clockwise direction black. 
The hexagons in the interior of $D_\delta$ 
are then randomly colored black or white 
with equal probability. This choice of boundary conditions forces there 
to be an interface which runs between $z_\delta$ and $w_\delta$. 
This interface, known as the percolation exploration process,  
has been proven to converge in distribution 
to the SLE$_6$ trace \cite{smirnov,camia_newman}.
The SLE$_6$ trace does not cross itself, but it does have self
intersections where the curve touches itself without crossing \cite{RS}. 
So the SLE$_6$ trace forms loops. 
Before the scaling limit, the percolation explorer does
not intersect itself, and so it does not form loops. 
However, it does often return to within one lattice spacing of itself and
so forms ``near-loops.''
The new stochastic process that we study is defined by erasing these 
near-loops in chronological order. 

Our loop-erased percolation explorer is similar in construction to the
loop-erased random walk (LERW), so we begin with a brief review of the
LERW \cite{Lawler_LERW}.
It can be defined on any lattice in any number of dimensions. 
Take a bounded domain containing the origin. We start a random walk on 
the lattice at the origin and stop the walk when it exits the domain. 
The walk can return to sites it has visited before and so form loops. 
We erase the loops it forms in chronological order. More precisely, 
if $\omega(0), \omega(1), \cdots, \omega(n)$ are the sites in the random 
walk, then its loop-erasure $\eta(0), \eta(1), \cdots, \eta(m)$ is 
defined as follows. Start by defining 
\bea
t_0 = \max \{ i : i \le n, \omega(i)=\omega(0) \}
\eea
and $\eta(0)=\omega(t_0)=\omega(0)$. 
Suppose we have defined $\eta(0), \eta(1), \cdots, \eta(j)$ and 
$t_0, t_1, \cdots, t_j$. If $\omega(t_j)=\omega(n)$ we stop and set $m=j$.  
Otherwise we let 
\bea
t_{j+1}= \max \{ i: t_j < i \le n, \omega(i)=\omega(t_j+1) \} 
\eea
The max in the above is the last time that $\omega$ visits $\omega(t_j+1)$. 
It is possible that this last time is just $t_j+1$, in which case $t_{j+1}$ 
is just $t_j+1$. Finally, we let $\eta(j+1)=\omega(t_{j+1})=\omega(t_j+1)$. 
In two dimensions the LERW has been proved to converge to 
radial SLE$_2$ in the scaling limit \cite{lerw_sle}.

\begin{figure}[tbh]
\includegraphics{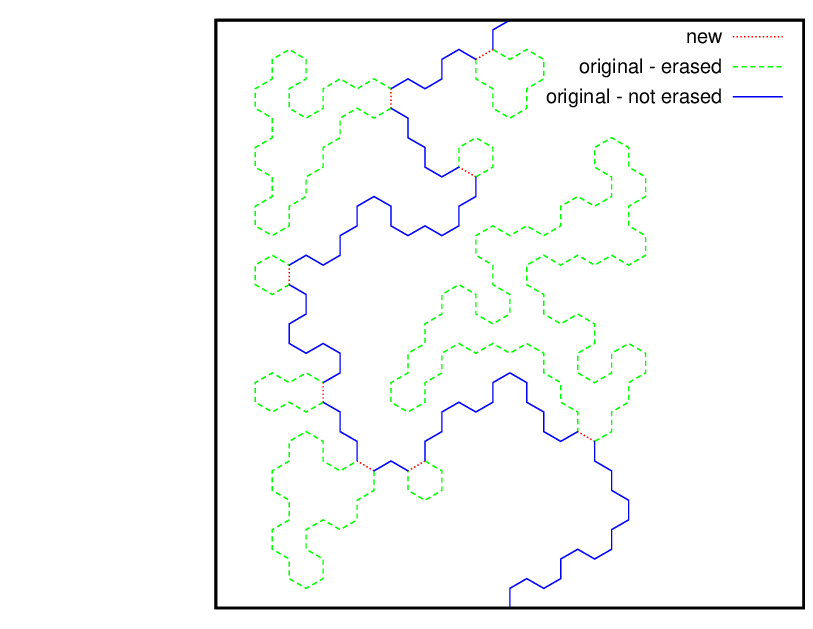}
\caption{The loop-erasure process. The original percolation interface is 
made up of the solid (blue) and dashed (green) bonds. The latter are erased
and replaced by the dotted (red) bonds. So the final curve is made up of the 
solid (blue) and dotted (red) bonds. 
}
\label{loop_erasing}
\end{figure}

\begin{figure}[tbh]
\includegraphics{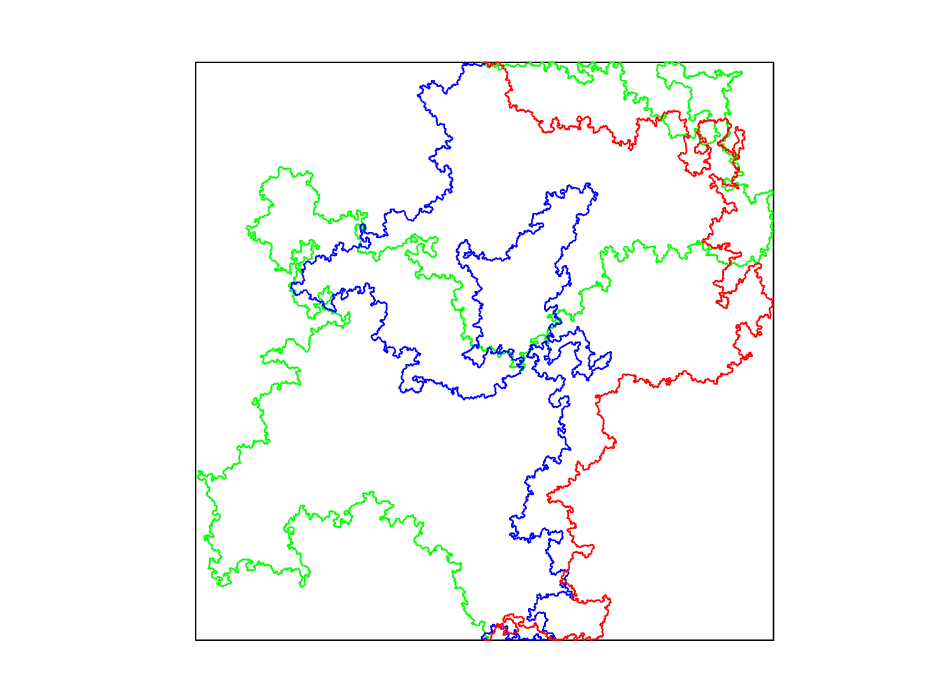}
\caption{Three samples of the loop-erased percolation explorer}
\label{samples}
\end{figure}

Now consider a percolation explorer path. When it returns to within
one lattice spacing of itself we erase this near-loop and replace it with
a single bond. 
More precisely, our loop-erasure process is defined inductively as follows. 
If $\omega(0), \omega(1), \cdots, \omega(n)$ are the sites in the percolation
explorer, then its loop-erasure $\eta(0), \eta(1), \cdots, \eta(m)$ is 
defined as follows. Start by defining $\eta(0)=\omega(0)$ and $t_0=0$. 
Suppose we have defined $\eta(0), \eta(1), \cdots, \eta(j)$ and 
$t_0, t_1, \cdots, t_j$. If $\omega(t_j)=\omega(n)$ we stop and set $m=j$.  
Otherwise we define
\bea
t_{j+1} = \max \{ i : t_j < i \le n, |\omega(i)-\omega(t_j)|=\delta \}
\eea
(Recall that $\delta$ is the lattice spacing.)
On the hexagonal lattice each site only has three nearest neighbors. 
Two of the nearest neighbors of $\omega(t_j)$ are occupied by 
$\omega(t_j-1)$ and $\omega(t_j+1)$. So the set we are taking the max of 
always contains $t_j+1$, and it can contain at most one other time. 
We then define $\eta(j+1)=\omega(t_{j+1})$. 
The loop-erasure process for the percolation explorer is illustrated in 
figure \ref{loop_erasing}. 
We will refer to the curve we get by this loop-erasure for the percolation
explorer as the loop-erased percolation explorer. 
Figure \ref{samples} shows a few samples of the loop-erased percolation 
explorer in a square. The lattice spacing is $1/1000$ of the length of the 
side of the square. 

The percolation explorer starting at $z$ and ending at $w$ generates 
the same path as the percolation explorer starting at $w$ and ending at $z$. 
However, it is easy to see from considering examples that the loop-erased 
percolation explorer depends on the direction in which we traverse the path.
We will always use $z$ to label the starting point and $w$ to label the 
terminal point. For the LERW simple examples show that we can get a different 
path when we loop-erase the random walk in reverse chronological order. 
Nonetheless, the distribution of the LERW using
reverse chronological order is the same as the distribution of the original 
LERW \cite{LERW_reversible}.
This property is known as reversibility. 
Our simulations indicate that the loop-erased percolation explorer is 
reversible in the scaling limit. We have not investigated if the discrete 
model is reversible before we take the scaling limit. 

When we erase the near-loops in the percolation explorer we obtain a new curve
which is smoother than the original percolation explorer in the 
sense that it has smaller fractal dimension.
(The SLE$_6$ trace has fractal dimension $7/4$ \cite{beffara}.)
The relation of the loop-erased percolation explorer to the original
percolation explorer is similar to the relation of the perimeter of a
percolation cluster to the external or accessible perimeter of the cluster
defined as follows. The cluster will have deep fjords which are connected to
the complement of the cluster only through an opening
whose width is on the order of a lattice spacing.
If we fill in these deep fjords, the perimeter of
the resulting object is called the external perimeter of the cluster.
It can also be defined by considering adsorbent particles with a diameter
that is slightly larger than the lattice spacing \cite{ga1,ga2,aizenman}.
If we follow the perimeter of a percolation cluster, 
the near-loops will be of two types - those forming a deep fjord
into the cluster and those forming a blob that is just barely attached
to the cluster. If we only erase the near-loops forming fjords, we will
obtain the external perimeter of the cluster.

\clearpage

\section{Tests of conformal invariance}

\begin{figure}[tbh]
\includegraphics{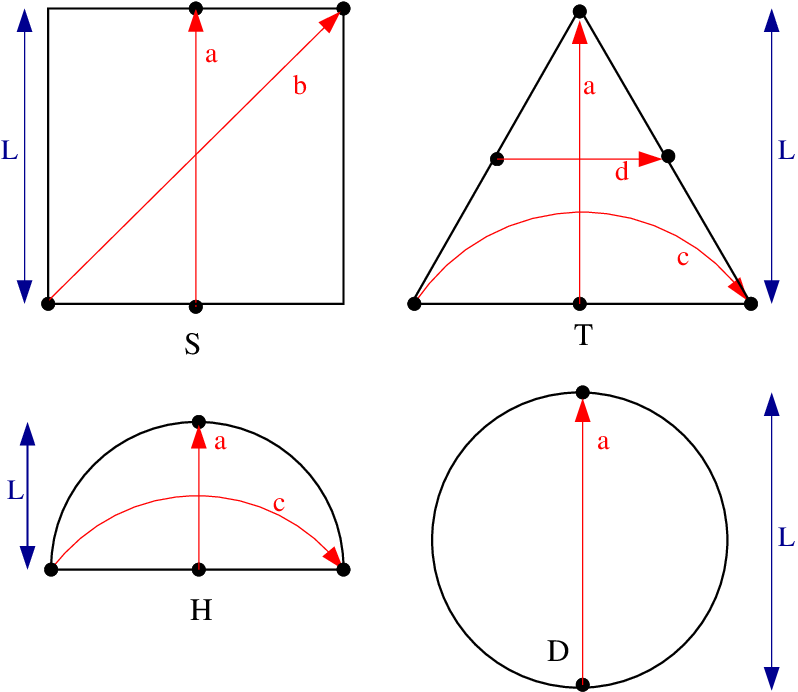}
\caption{The four regions we study: S=square, T=equilateral triangle,
  H= half-disc, D=disc. The red arrows indicate difference choices of the
  starting and terminal points for the explorer.
  }
\label{regions}
\end{figure}

We test the conformal invariance of the loop-erased percolation explorer 
by simulating it for four different domains
(square, equilateral triangle, half-disc and disc) which we will denote by 
S,T,H,and D. We also use 
several choices of the starting point $z$ and terminal point $w$. 
They are shown in figure \ref{regions}.
For the square we have two choices of  starting and terminal points. 
They are labeled a and b in the figure. 
For the triangle T we have four choices of starting and terminal points. 
Choices a,c and d are shown in the figure. Choice b is given by reversing the 
direction of choice a. For the half-disc H, choices a and c are shown in the 
figure. Choice b is the reversal of a. Finally, for the disc D there is only one
choice of starting and terminal points. 
So there are ten choices of region and starting/terminal points which we
label Sa, Sb, Ta, Tb, Tc, Td, Ha, Hb, Hc, Da. 

We can take the scaling limit by fixing the domain and introducing a 
lattice with spacing $\delta$ which is then sent to zero, 
or by taking the lattice spacing to be $1$ and 
letting the scale of the domain go to infinity. We do the latter. 
The length $L$ is indicated for each domain in the figure.  
For the triangle, half-disc and disc we have done simulations with 
$L=200, 283, 400, 566, 800$. These values are chosen so that $L$ is increasing
by approximately a factor of $\sqrt{2}$. 
The domain in the simulation is made up of hexagons, so it is only an 
approximation to the true domain. As $L$ varies the approximating domain can 
change in a somewhat irregular way. As a result the finite $L$ effects 
can show a somewhat chaotic dependence on $L$. As we will see in detail later,
this chaotic dependence is not that significant for the triangle, half-disc
and disc. However, it can be quite pronounced for the square. The reason is 
that the approximating domain is a rectangle of hexagons, but the aspect 
ratio of this rectangle is not $1$. As $L$ varies this aspect ratio changes 
enough to produce noticeable chaotic dependence of the finite $L$ 
effect on $L$. 
We have found that this can be greatly reduced by choosing values of 
$L$ for which the aspect ratio does not vary so much. The values 
$L=204, 273, 405, 564, 810$ do this, so these are the values we use
for the simulations for the square.

We test conformal invariance in two ways. 
The first uses a family of random variables. 
Fix a conformal map $\phi$ from the domain to the upper half plane which 
sends the starting point $z$ to the origin and the terminal point $w$ to 
$\infty$. (There is a one parameter family of such maps.)
We let $C_r$ be the curve in the domain whose image under
$\phi$ is a semicircle of radius $r$ centered at the origin. 
We find the first point $p$ where the loop-erased percolation 
explorer crosses the curve $C_r$. The random variable is the polar angle of 
$\phi(p)$. (For convenience we divide this angle by $\pi$.) 
Just how this random variable depends on $r$ depends on the choice of $\phi$,
so we do not parameterize the random variable by $r$ but rather by a parameter
$t$ defined as follows. Consider the intersection of $C_r$ with the line
from the starting point to the terminal point. Then $t$ is the distance from
the starting point to this intersection divided by the distance from the
starting point to the terminal point.
If the scaling limit is conformally invariant,
then the distribution  of this random variable will be the same for 
all domains, all choices of $z,w$ and all choices of the parameter $t$.  
We will refer to this random variable as the ``first hit'' random variable. 

The second test of conformal invariance uses 
the probability of passing right of a point $p$ in the 
domain. If we take a conformal map $\phi$ from our domain to the upper 
half plane which sends $z$ to $0$ and $w$ to $\infty$, then this probability
only depends on the polar angle $\theta$ of $\phi(p)$. 
Rather than compute this probability for a single point, we compute it 
for a sequence of points along a line segment in the domain.
We parameterize the line segment by $\theta$, 
and look at the probability of passing 
right of the points on the line as a function of $\theta$. 
We will refer to this function as the ``pass right function.''
If the model is conformally invariant then 
the pass right function will be the same for 
all domains, all choices of $z,w$ and all choices of the line in the domain. 
Note that the image of the line under the conformal map
is usually not a semi-circle in the half-plane, but 
this does not matter since the probability of passing right of a point 
in the upper half plane only depends on the polar angle, not on the radius. 
The line segments we use for our different choices of domains and 
starting points are all horizontal or vertical. For
Sa, Sb, Ta, Tb, Ha, Hb, and Da,  
the line segment is horizontal. For Tc, Td and Hc it is vertical. 
The position of the line segment is parameterized by $t$. The position is 
linear in $t$ with $t=0$ corresponding to the line segment passing through
the point where the explorer starts and $t=1$ to the line segment
passing through the termination point.

Before we take the scaling limit, the random variable we are studying is 
discrete since there are only a finite number of points where the
loop-erased percolation explorer can first hit the curve $C_r$.
So the cumulative distribution function (cdf)
is a step function. Similarly, the 
pass right function is a step function before the scaling limit. 
In the scaling limit these step functions should converge to smooth functions,
but for the lattice spacings that can be simulated the effect of this 
discreteness is quite noticable.
We can reduce the effect of this discreteness in the following way. 
Rather than consider the first hit random variable for a single value of 
$t$, we average the random variable over some interval for $t$. 
If the loop-erased percolation explorer is conformally invariant, then 
the cdf of this averaged random variable will be independent of the interval 
we average over, as well as the domain and starting and terminal points. 
As $t$ varies the finite set of $\theta$ where the cdf 
jumps changes, so this averaging over $t$ 
reduces the effect of the discreteness of the random variable for a fixed $t$.
Similarly we can reduce the effect of the discreteness for the pass right 
function by averaging $t$ over an interval. 
For both the first hit random variable and the pass right function
we average $t$ over three 
intervals: $[0.25,0.35]$, $[0.45,0.55]$ and $[0.65,0.75]$.

\begin{figure}[tbh]
\includegraphics{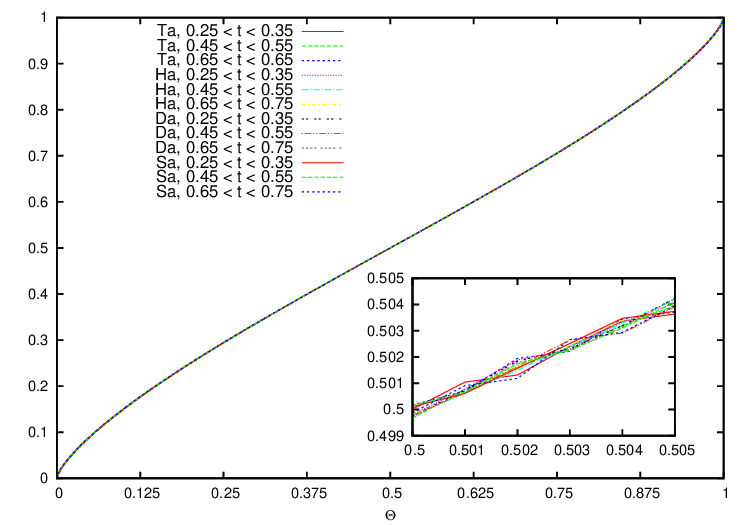}
\caption{The cdf for the first hit random variable for the four domains and
three choices of the interval for $t$. The 12 curves are indistinguishable 
in the main figure. The inset blows up a tiny portion of the main plot to see
their difference.}
\label{cdf_all}
\end{figure}

\begin{figure}[tbh]
\includegraphics{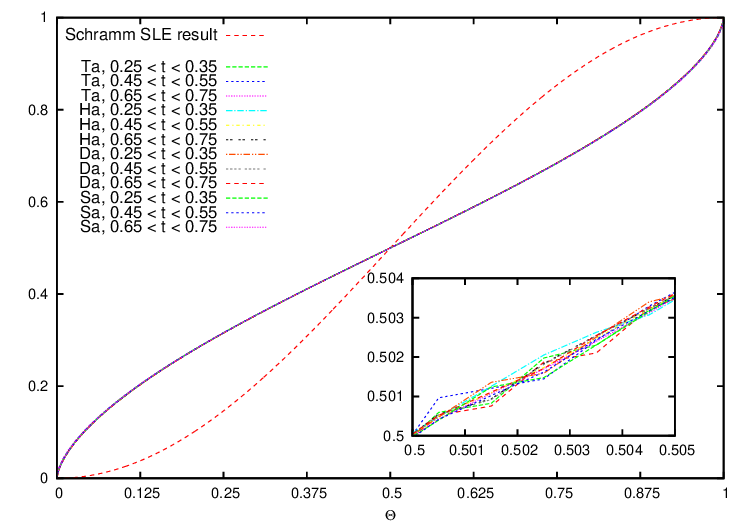}
\caption{The pass right function for the four domains and
three choices of the interval for $t$. The inset blows up a
tiny portion of the main plot to see the difference in the 12 curves.
The dashed curve is the exact function for SLE$_{8/3}$. 
}
\label{pfunc_all}
\end{figure}

\begin{figure}[tbh]
\includegraphics{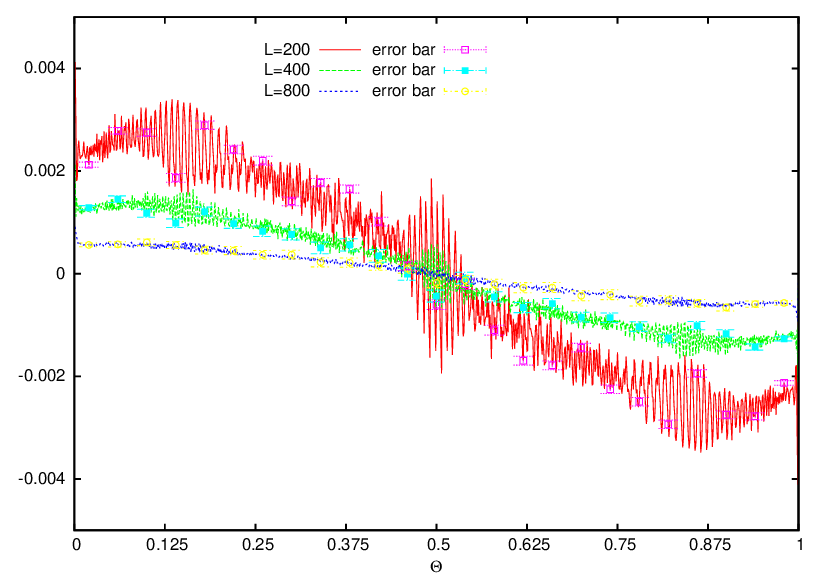}
\caption{
The differences between the cdf's for the first hit random variable
and the ad hoc reference function for Ta with $L=200, 400, 800$.}
\label{40_dif2}
\end{figure}

In our simulations to test the conformal
invariance we generated $10^8$ samples for each of the thirty cases
(10 possibilities for the region and starting/terminal points and
3 possibilities for the interval over which we average $t$).
When we compute a cdf or the pass right function we are computing a
probability for each value of $\theta$. Since our samples are independent,
the variance for our estimate is $p(1-p)/N$. So for $p$ around $1/2$,
two standard deviations is approximately $10^{-4}$. We have not included
these error bars in our plots of the first hit random variable cdf or the
pass right function to keep the figures from being too cluttered,
except in figure \ref{40_dif2}.

Figure \ref{cdf_all} shows the cdf of the 
first hit random variable for Sa, Ta, Ha and Da and for all three
different intervals for $t$. 
There are $12$ curves in the figure, but they look identical
in the main plot.  The inset blows up a tiny portion of the main plot
to illustrate the size of the differences in the $12$ curves.
In the inset the differences are roughly $1/1000$, and this is typical for
all $\theta$.
Figure \ref{pfunc_all} shows the pass right function for the same four cases
of regions and starting/terminal points
and all three choices of intervals for $t$.  
Again, the $12$ curves in the main plot are indistinguishable.
The differences for these 12 curves are also roughly $1/1000$. 
The dashed curve in the figure 
is the exact result for the pass right function for SLE$_{8/3}$.  

\begin{figure}[tbh]
\includegraphics{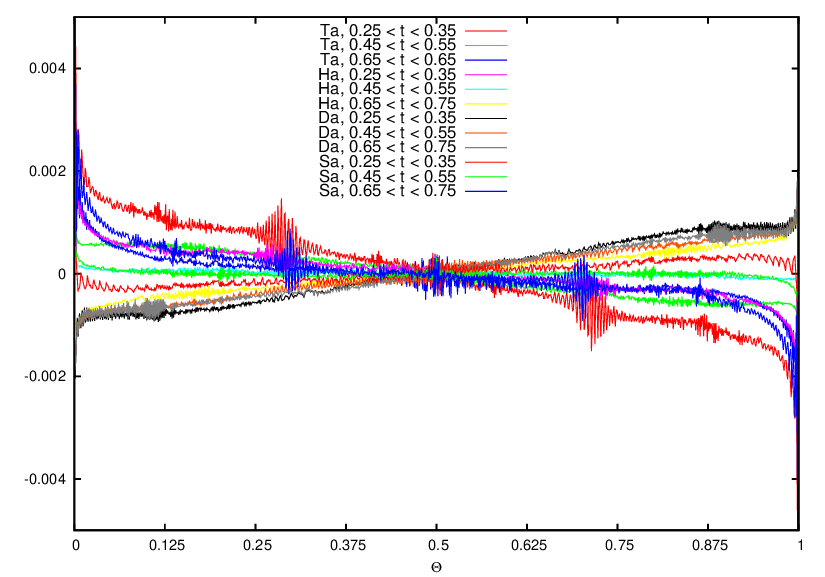}
\caption{
We plot the difference between the 12 first hit random variable
cdf's plotted in figure \ref{cdf_all}
and an ad hoc reference function. 
}
\label{dif_all}
\end{figure}

\begin{figure}[tbh]
\includegraphics{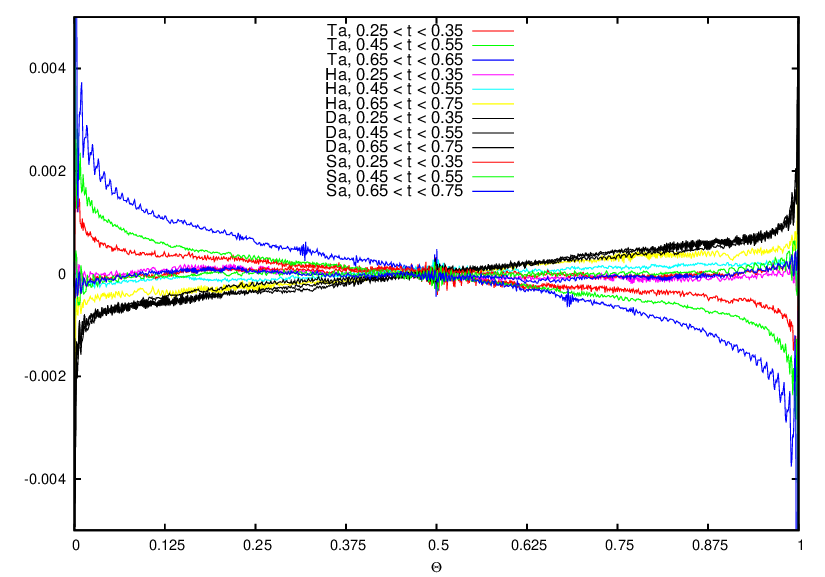}
\caption{
We plot the difference between the 12 pass right functions
plotted in figure \ref{pfunc_all} and an ad hoc reference function. 
}
\label{pdif_all}
\end{figure}

To ``zoom in'' on these plots we will subtract off a reference function. 
We do not have conjectures for what this first hit cdf and the pass right 
function are, so for our reference functions we just average four of the
cases.  We use Sa, Ta, Ha, Da with the
parameter $t$ averaged over $[0.45,0.55]$.
Figure \ref{dif_all} shows the cdf minus the reference function for the ten
choices of region and starting/terminal points. 
We only show the curves for the simulations with the
parameter $t$ averaged over $[0.45,0.55]$.
Figure \ref{pdif_all} shows the pass right functions for the same ten cases
minus the pass right reference function.

The results in figures \ref{dif_all} and \ref{pdif_all} only
use a single value of $L$ and there is no attempt to extrapolate to
$L=\infty$. Figure \ref{40_dif2} shows the difference between the cdf
for $L=200, 400, 800$ and the reference function for Ta. 
Even with the averaging of the parameter $t$, the cdf's and
pass right functions are not smooth enough to extrapolate these functions
point wise. We will instead study their Fourier coefficients for several
values of $L$ and attempt to extrapolate them to $L=\infty$.

\clearpage

\clearpage

Our functions are defined on $[0,1]$ and we compute the Fourier
coefficients in the expansion
\bea
a_0 + \sum_{n=1}^\infty a_n \cos(2 \pi n x)
+ \sum_{n=1}^\infty b_n \sin(2 \pi n x)
\eea
For the RV we compute the Fourier coefficients of the density $\rho(x)$
rather than the cdf.
Since the density has the symmetry $\rho(1-x)=\rho(x)$, the $b_n$ are all zero.
And since it has integral $1$, $a_0=1$. So we only compute $a_1, a_2, \cdots$.
The pass right function satisfies $p(1-x)=1-p(x)$. This implies $a_0=1/2$
and all the rest of the $a_n$ are zero. So we only compute
$b_1,b_2,\cdots $. 

\begin{figure}[tbh]
\includegraphics{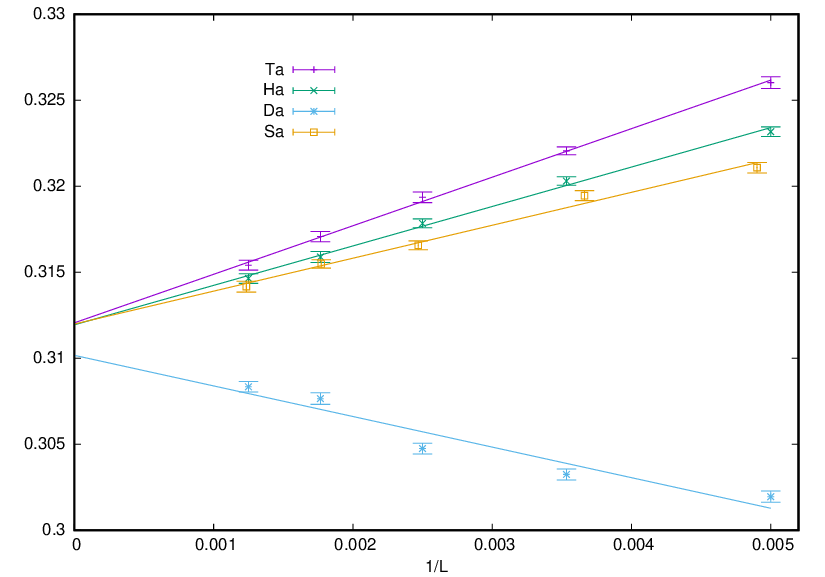} 
\caption{
The Fourier coefficient $a_1$ for the first hit random variable cdf
as a function of $1/L$. Four cases (Ta, Ha, Da, Sa) are plotted.
The lines are least squares fits to the data. 
}
\label{fc_rv_1}
\end{figure}

\begin{figure}[tbh]
\includegraphics{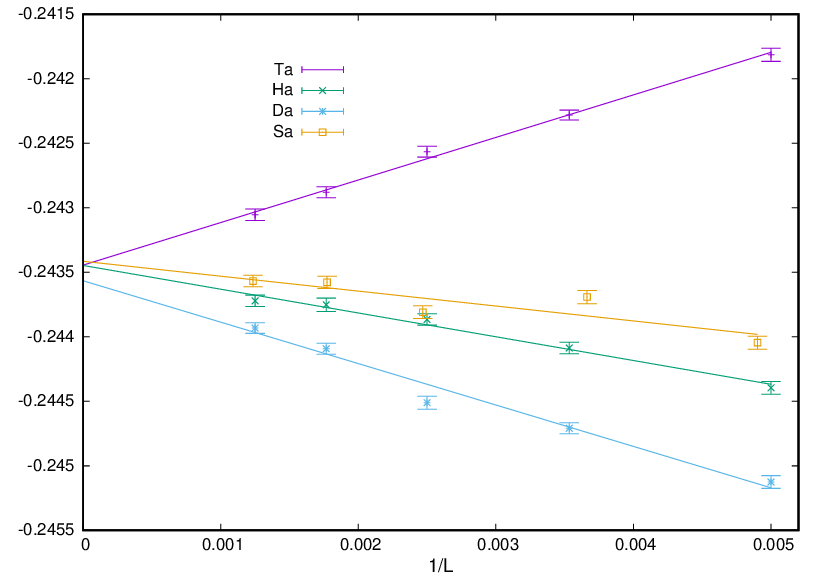}
\caption{
The Fourier coefficient $b_1$ for the pass right function
as a function of $1/L$. Four cases (Ta, Ha, Da, Sa) are plotted.
The lines are least squares fits to the data. 
}
\label{fs_pfunc_1}
\end{figure}

As $L$ changes, the way in which our domains are approximated
with hexagons changes in
a somewhat erratic way. So even the Fourier coefficients have a somewhat
chaotic dependence on $L$. In figure \ref{fc_rv_1} we plot the Fourier
coefficient $a_1$ for the cdf as a function of $1/L$ for regions
Ta, Ha, Da, Sa with the parameter $t$ averaged over
$[0.45,0.55]$. The lines are least squares fits to the data. 
We will use the intercept with the vertical axis as the extrapolation of
the coefficient to $L=\infty$. 
Figure \ref{fs_pfunc_1} shows the analogous plots for the Fourier
coefficient $b_1$ for the pass right function.

\clearpage

\begin{figure}[tbh]
\includegraphics{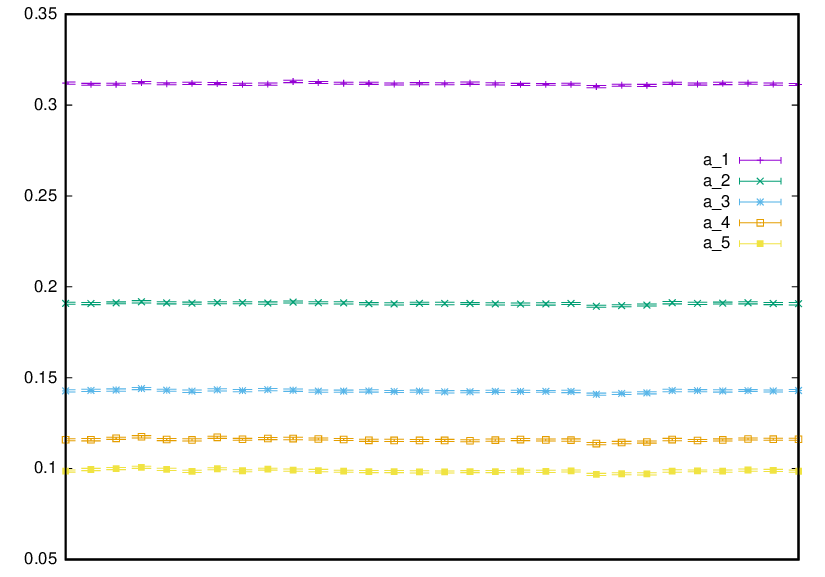}
\caption{
For the density of the first hit random variable
we plot the five largest Fourier coefficients 
for all ten choices of domain and $z,w$ and all three choices of the
interval for averaging $t$.
}
\label{fc_rv_intercept_all}
\end{figure}

\begin{figure}[tbh]
\includegraphics{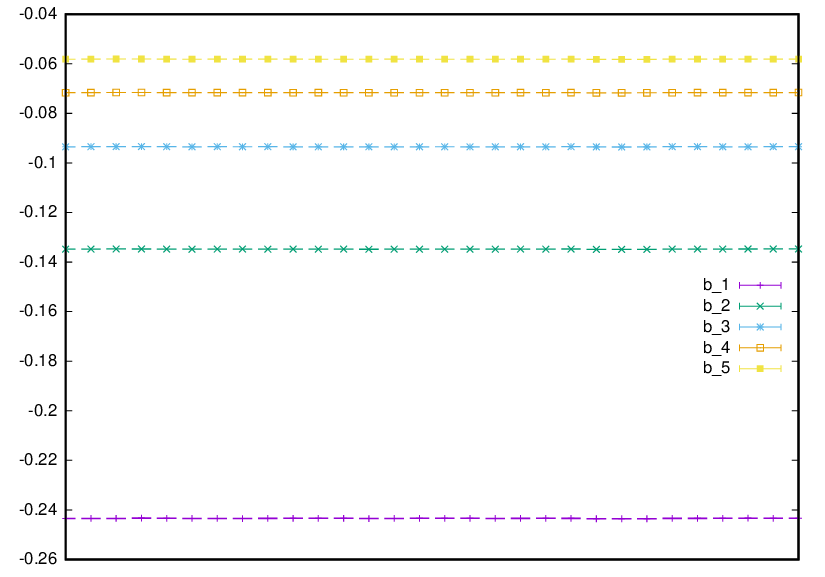}
\caption{
For the pass right function 
we plot the five largest Fourier coefficients 
for all ten choices of domain and $z,w$ and all three choices of the
interval for averaging $t$.
}
\label{fs_pass_intercept_all}
\end{figure}

Our next two plots show the extrapolated Fourier coefficients 
for all ten choices of domain and $z,w$ and all three choices of the
interval for averaging $t$ for  
the first hit random variable density (figure \ref{fc_rv_intercept_all})
and the pass right function (figure \ref{fs_pass_intercept_all}).
Going from left to right, the order of the domains is 
Ta, Tb, Tc, Td, Ha, Hb, Hc, Da, Sa, Sb with three data points for each domain
corresponding to averaging the parameter $t$ over the usual three intervals 
- $[0.25,0.35]$, $[0.45,0.55]$ and $[0.65,0.75]$.
We plot the extrapolated values for the five largest Fourier coefficients
for these thirty cases. If the model is conformally invariant then all
thirty cases should give the same values for the Fourier coefficients.
As the plots show, the values for the Fourier coefficients are very nearly
the same. The error bars shown in the plots are the statistical errors
from the Monte Carlo simulation. There is also error arising from the
extrapolation to $L=\infty$. This is difficult to estimate given the somewhat
chaotic dependence of the Fourier coefficients on $L$. 

\clearpage

\section{Dimension of the curve}

\begin{figure}[tbh] 
\includegraphics{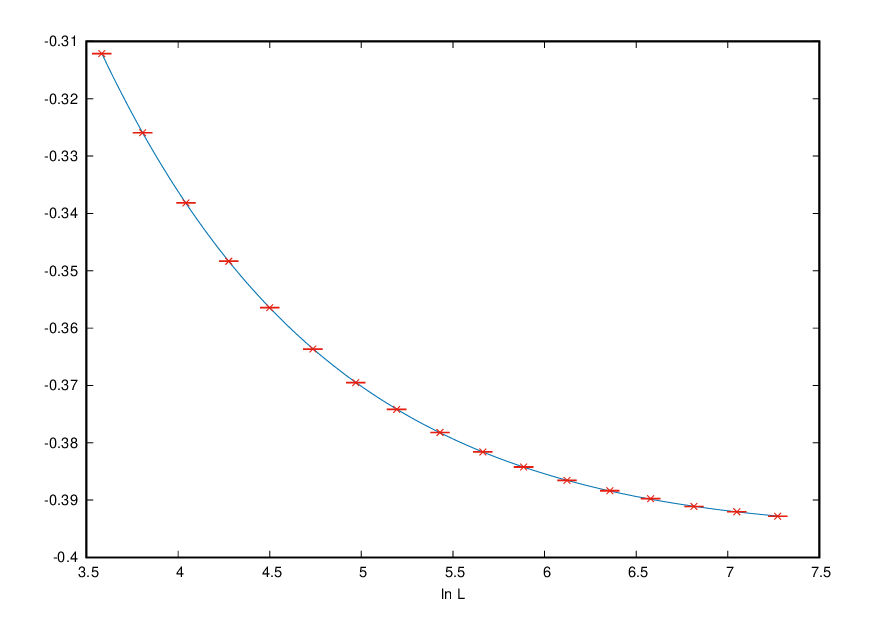}
\caption{Log-log plot of average number of steps to travel distance $L$ with
  $4/3 \ln(L)$ subtracted off.}
\label{nu_fit}
\end{figure}

The average distance the loop-erased percolation explorer
travels as a function of the number of steps $N$
should be asymptotically proportional to $N^\nu$ for some exponent $\nu$.
The fractal dimension of the curve is $1/\nu$. 
In our simulations the distance between the starting and terminal points
of the loop-erased percolation explorer is fixed and
the number of steps it takes to travel that 
distance is random. The average number of steps 
should be asymptotically proportional to  $L^{1/\nu}$.
We have computed this average number of steps for the triangular region
for $17$ values of $L$ ranging from $36$ to $1440$ with $5 \times 10^7$ samples
for each value of $L$. 

To estimate $\nu$ accurately we must take into account the next order
term. Let $\overline{N}(L)$ be the average number of steps. We assume 
$\overline{N}(L) = c L^{1/\nu} (1 + a L^{-\delta} + \cdots)$.
So 
\bea
\ln(\overline{N}(L)) &=& \ln(c) + \frac{1}{\nu} \ln(L) +\ln(1 + a L^{-\delta}) 
+ \cdots \label{nu_fit_eq}
\nonumber
\\
 &=& \ln(c) + \frac{1}{\nu} \ln(L) + a L^{-\delta} + \cdots
\eea
This is linear in the unknown parameters 
$\ln(c), \frac{1}{\nu}$ and $a$, but not in $\delta$. 
For a given value of $\delta$ we do a weighted least squares fit to 
find $\ln(c), \frac{1}{\nu}$ and $a$. We then search over $\delta$ 
to find the value that minimizes the residual sum of squares (RSS). 
We find the RSS is minimized when $\delta=0.7552$ and 
for this $\delta$, $1/\nu = 1.334782 \pm 0.000038$.
We emphasize that the error bars on $1/\nu$ are only the error from the 
Monte Carlo. There is also error from the neglected higher order terms in eq. 
\reff{nu_fit_eq}.
Figure \ref{nu_fit} shows a plot of $\ln(\overline{N}(L)) - \frac{4}{3} \ln(L)$
as a function of $\ln(L)$. The curve shown is 
$ \ln(c) + (\frac{1}{\nu}-\frac{4}{3}) \ln(L) + a L^{-\delta}$. It fits the
data quite well, supporting our ansatz \reff{nu_fit_eq}.

Since the numerical estimate of $\nu$ is very close to $3/4$, if 
the loop-erased percolation explorer is some SLE$_\kappa$, then $\kappa$ must
be very close to $8/3$. For SLE, Schramm found an explicit formula 
for the pass right function \cite{schramm_pass_right}. Figure \ref{pfunc_all}
includes Schramm's result for $\kappa=8/3$ (the dashed curve). 
It is clearly different from the pass right function computed in our 
simulations, leading us to conclude that the loop-erased 
percolation explorer is not an $SLE$. 

\section{Loop-erasing with a gap of $2$}

Our definition of the loop-erasure process was that when the percolation
explorer comes within one lattice spacing of itself, we erase this near
loop and replace it with a single bond. From now on we will refer to this
as loop-erasure with a gap of $1$. Now we will consider loop-erasure with
a gap of $2$. Loosely speaking we consider the percolation
explorer to have formed a loop when it comes within two lattice spacings of
itself. We erase such a near loop and replace it with two bonds.
(Note that unlike the square lattice, on the hexagonal lattice there is only
one choice for the two bonds.)

However, the above definition is a bit too simplistic. Suppose that
$s<t$ are times such that $|\omega(t)-\omega(s)|=2 \delta$, 
i.e., the loop at time $t$ has returned to within $2$ lattice spacings 
of its past.
Suppose that the next step after time $t$ brings the loop within
$1$ lattice spacing of $\omega(s)$, i.e.,  $|\omega(t+1)-\omega(s)|=\delta$.
Then replacing the near-loop from $\omega(s)$ to $\omega(t)$ with the
two bonds from $\omega(s)$ to $\omega(t)$ will result in a path in which
the bond from $\omega(t)$ to $\omega(t+1)$ is traversed twice.
So we define loop-erasure with a gap of $2$ as follows.

We start by defining $\eta(0)=\omega(0)$.
Unlike the definition when the gap is $1$, it is no longer true that all 
the sites in the loop-erased explorer are sites in the original explorer.
For $i$ such that $\eta(i)$ is a site in $\omega$, we define $t_i$ by
$\omega(t_i)=\eta(i)$. 
Now suppose we have defined $\eta(0), \eta(1), \cdots, \eta(j)$, and $j$ is
such that $\eta(j)$ is a site in $\omega$.  
If $\eta(j)=\omega(n)$ we stop and set $m=j$.  
Note that $\eta(j-1)$ and $\omega(t_j+1)$ are two of the three
nearest neighbors of 
$\omega(t_j)$. Let $x$ denote the third nearest neighbor of $\omega(t_j)$. 
We first check if there is a $k>t_j+1$ such that 
$|\omega(k)-x|=\delta$, i.e., a gap of size $1$. 
If so, we set $\eta(j+1)=\omega(k)$ and $t_{j+1}=k$. 
In this case the definition of the loop-erased explorer is extended by one 
step.
If there is not such a $k$, we check if there is a $k>t_j+1$ such that 
$\omega(k)$ is a nearest neighbor of $x$. If there is, then 
we set $\eta(j+1)=x$ and $\eta(j+2)=\omega(k)$. 
We have $t_{j+2}=k$ and $t_{j+1}$ is not defined since $\eta(j+1)$ is 
not a site in $\omega$. In this case the definition of the loop-erased 
explorer is extended by two steps.
Finally, if we did not find a gap of size $1$ or $2$, then we set 
$\eta(j+1)=\omega(t_j+1)$ and so $t_{j+1}=t_j+1$. 

\begin{figure}[tbh] 
\includegraphics{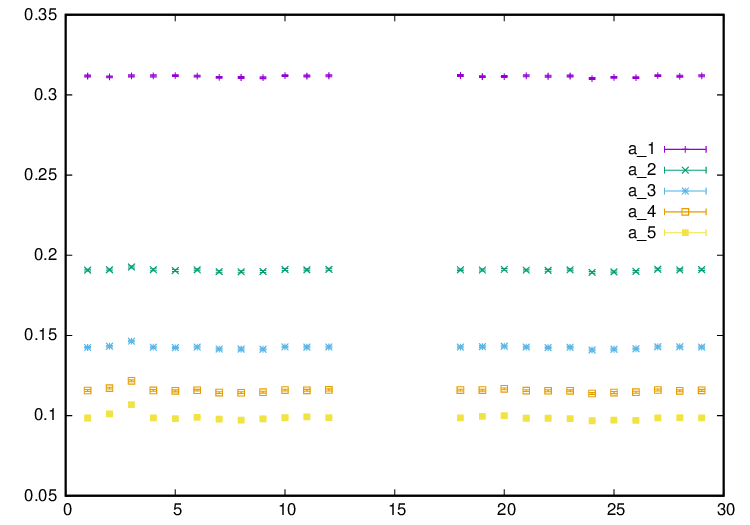}
\caption{
For the density of the first hit random variable 
we plot the five largest Fourier coefficients 
for four choices of domain and $z,w$ and all three choices of the
interval for averaging $t$. The data on the left is for loop-erasing with
a gap of $2$, and on the right for loop-erasing with a gap of $1$. 
}
\label{gap2_rv_fc}
\end{figure}

\begin{figure}[tbh] 
\includegraphics{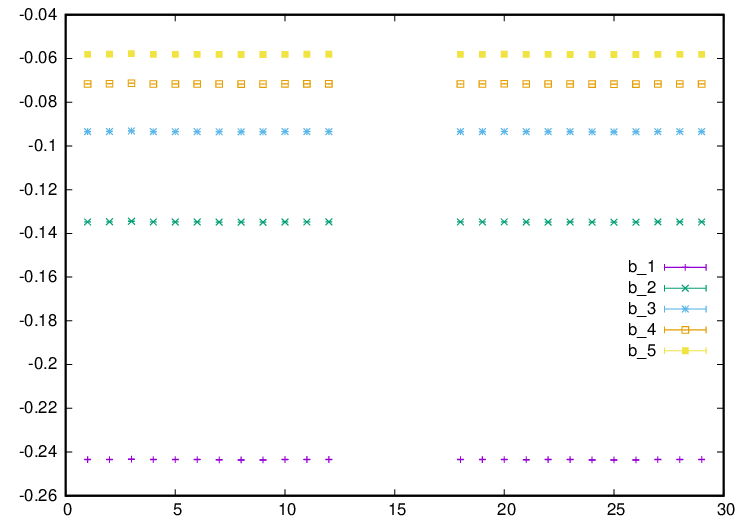}
\caption{
For the pass right function 
we plot the five largest Fourier coefficients 
for four choices of domain and $z,w$ and all three choices of the
interval for averaging $t$. The data on the left is for loop-erasing with
a gap of $2$, and on the right for loop-erasing with a gap of $1$. 
}
\label{gap2_pass_fs}
\end{figure}

In this section we study whether the loop-erased percolation explorer using
a gap of $2$ has the same scaling limit as  the loop-erased
percolation explorer using a gap of $1$.
We study this question by computing the same two quantities that we
did for the case of a gap of $1$, namely, the 
first hit random variable and the pass right function.
We compute the Fourier coefficients of the density of the random variable
and of the pass right function for the same set of lengths as we did
for a gap of 1. Then we extrapolate them to $L=\infty$ just as we did
before. We have only done these computations for 
Sa, Ta, Ha and Da and for all three intervals for averaging $t$.
So we have 12 cases instead of 30 as before. 
In figure \ref{gap2_rv_fc} we plot the extrapolated values for the
five largest Fourier coefficients for the density of the first hit
random variable. The 12 cases for erasing using a gap of $2$ are on the left.
From left to right the domains are Ta, Ha, Da, Sa with three data points for 
each domain for the three interval for $t$. 
On the right we show our previous results for erasing using a gap of $1$
for the same 12 cases.
The analogous plot for the five largest Fourier coefficients for the
pass right function is figure \ref{gap2_pass_fs}.
As can be seen in these two figures, the results for a gap of 1 and a
gap of 2 appear to be the same.

\clearpage

\section{Macroscopic difference of loop-erasures}

Suppose we take a percolation exploration process curve $\gamma_L$
and loop-erase it
in two different ways - one using a gap of $1$ and the other a gap of $2$.
We let $\gamma^1_L$ and $\gamma^2_L$ be the resulting two curves.
As before we take the lattice spacing to be $1$ and take the scaling limit
by letting $L \rightarrow \infty$. So we need to rescale the curves by
a factor of $1/L$ to take the scaling limit. 
Note that for a given $L$, $\gamma^1_L$ and
$\gamma^2_L$ are random curves on the same probability space.
The conclusion of the previous section is that $L^{-1} \gamma^1_L$ and
$L^{-1} \gamma^2_L$ converge in distribution to the same process.
So the distributions of $L^{-1} \gamma^1_L$ and $L^{-1} \gamma^2_L$
are close when $L$ is large. 
Since these two curves are defined on the same probability space (which
depends on $L$), we can ask if the two curves are close with high probability.
(The convergence in distribution does not imply that they must be close.)
More precisely we can ask if for all $\epsilon>0$ we have 
\bea
\lim_{L \rightarrow \infty} P(d(L^{-1} \gamma^1_L,L^{-1} \gamma^2_L) \ge \epsilon) = 0
\label{cauchy_in_prob}
\eea
where $d()$ is some distance function for curves.

Computing the distance between two curves is computationally intensive
since we must consider all possible parameterizations of the curves. 
So we will study a different quantity to test if $\gamma^1_L$ and
$\gamma^2_L$ are close. 
Let $C$ be a curve which goes from one boundary point of the domain
to another boundary point in such a way that it disconnects the starting 
and ending points.  So the curves $\gamma^1_L$ and $\gamma^2_L$
must cross $C$ at least once. 
Let $t$ parameterize $C$ by arc-length, normalized so that $t$ runs
from $0$ to $1$. 
For the point corresponding to parameter value $t$, 
let $r(\gamma,t)$ be $1$ if $\gamma$ passes right of the point with parameter
value $t$, $0$ if it passes left.
The quantity $|r(\gamma^1_L,t)-r(\gamma^2_L,t)|$
is then the indicator function of the event that one of $\gamma^1_L$
and $\gamma^2_L$ passes right of the point and the other passes left of
the point. 
We will study the following random variable:
\bea
X_L = \int_0^1 |r(\gamma^1_L,t)-r(\gamma^2_L,t)| \, dt
\eea
So $X_L$ computes the fraction of the curve $C$ where one of 
$\gamma^1_L$ or $\gamma^2_L$ is right of the point and the other is 
left of the point. 

Our simulations will show that $E X_L$ does not converge to zero.
Strictly speaking, $X_L$ not being small does not imply that 
$d(L^{-1} \gamma^1_L,L^{-1} \gamma^2_L)$ is not small.
It is possible that the two curves
stay very close to each other but oscillate back and forth across $C$
in such a way that $X_L$ is not small, but this sort of behavior is not
expected. 
We will see in the simulations that for finite $L$,
$E X_L$ is nonzero but it is decreasing as $L \rightarrow \infty$. 
The tricky question is whether it converges to zero or not as
$L \rightarrow \infty$. 

We study $X_L$ just for Ta. We take the curve $C$ to be a horizontal line
which is a distance $tL$ above the base of the triangle. (Recall that $L$
is the height of the triangle.) The parameter $t$ is averaged over
$[0.4,0.6]$.
The equilateral triangle T has the advantage that it can be approximated by
hexagons in a regular way.
We have simulated a large number of values of $L$ to study the
$L \rightarrow \infty$ limit carefully.
We use $L=36$, $45$, $57$, $72$, $90$, $114$, $144$, $180$, $228$,
$288$, $360$, $456$, and $576$
with $5 \times 10^8$ samples computed for each value.
These lengths were chosen so that the ratio of consecutive lengths is very
nearly $2^{1/3}$. 

\begin{figure}[tbh]
\includegraphics{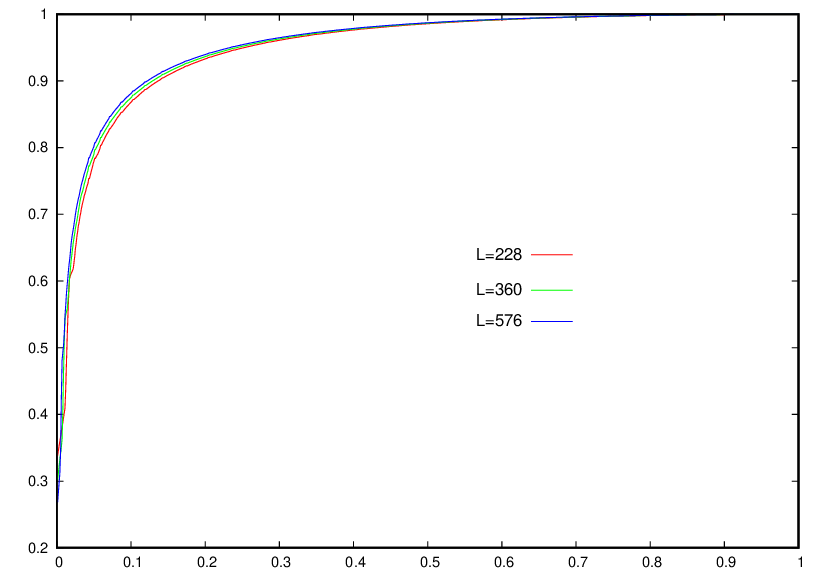}
\caption{The cdf of $X_L$ for $L=228, 360$ and $576$.}
\label{cdf}
\end{figure}

Figure \ref{cdf} shows the cdf of $X_L$ for $L=228, 360$ and $576$.
Since $E X_L$ is equal to $\int_0^\infty [1-F_L(t)] dt$ where $F_L(t)$ is the cdf
of $X_L$, the expected value $E X_L$ is equal to the area 
between the cdf of $X_L$, the horizontal line with
height $1$ and the vertical axis. We need to determine if this area goes to
$0$ as $L \rightarrow \infty$. 

\begin{figure}[tbh]
\includegraphics{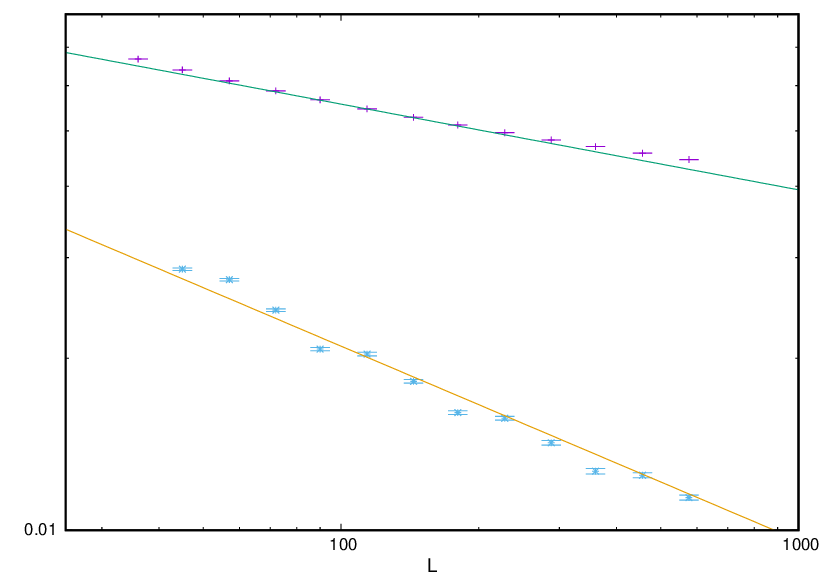}
\caption{The top curve is a log-log plot of $E X_L$ as a function of $L$.
The bottom curve is a log-log plot of the differences in $E X_L$ for
successive values of $L$. 
}
\label{means_both}
\end{figure}

If $E X_L$ goes to zero, it is natural to expect it goes as $L^{-p}$ for
some power $p$. So a log-log plot of $E X_L$ as a function of $L$ would
be linear. The top curve in figure \ref{means_both} shows this log-log plot.
It is not linear but rather is slightly convex.
The line shown as a guide to the eye has slope $-0.15$.
If $E X_L$ goes as $\mu + c L^{-p}$, then since our values
of $L$ are essentially growing geometrically, 
the differences between values of $E X_L$ for successive values of $L$
should go to zero as $L^{-p}$.
The lower curve in figure \ref{means_both} is a log-log plot of the
differences.
It is noisy but looks to be linear with slope around $-1/3$. Note that this
slope is quite different than the slope of the top line in the figure. 
We have multiplied the data for the bottom curve by a factor of $10$ so
that the two curves can be shown in a single figure.

\newpage

\begin{figure}[tbh]
\includegraphics{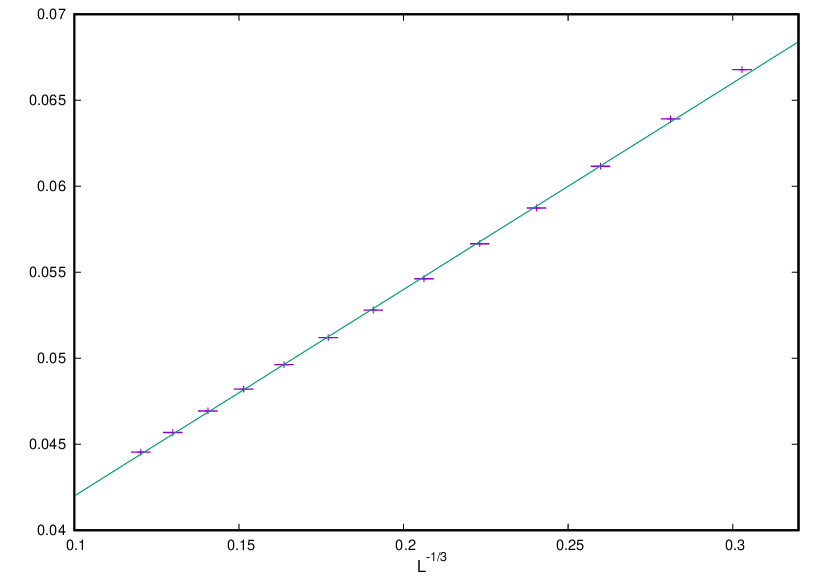}
\caption{
$E X_L$ is plotted as a function of $L^{-1/3}$. 
The nonzero vertical intercept shows that $E X_L$ does not converge to
zero as $L \ra \infty$. 
}
\label{meanslp}
\end{figure}

If $E X_L$  goes as $\mu + c L^{-p}$, then it should be a linear
function of $L^{-p}$. So we plot $E X_L$ as a function of $L^{-1/3}$
in figure \ref{meanslp}. Note that this is not a log-log plot.
The data is very well fit by a linear function with a vertical intercept that
is clearly not zero. This plot is the best evidence that $E X_L$  does
not converge to zero as $L \ra \infty$.

\clearpage

\section{Conclusions}

The bulk of our Monte Carlo simulations have been for the loop-erased
percolation explorer in which we consider the process to have formed a loop
when it comes within one lattice spacing of a site it has visited before
(the gap=1 model). These simulations give strong support to three conjectures:

\begin{conjecture} The scaling limit of the 
loop-erased percolation explorer is conformally invariant. 
\end{conjecture} 

\begin{conjecture} The fractal dimension of the scaling limit of the 
loop-erased percolation explorer is $4/3$. 
(This is the dimension of SLE$_{8/3}$.)
\end{conjecture} 

\begin{conjecture} 
The scaling limit of the loop-erased percolation explorer 
is \textbf{not} SLE$_{8/3}$. 
\end{conjecture} 

We have also carried out simulations in which the explorer is considered to
have formed a loop when it comes within two lattice spacings of a previously
visited site (the gap=2 model). These simulations support the conjecture
that the gap=1 and gap=2 models converge in distribution to the
same limit as the lattice spacing goes to zero. 
The gap=1 and gap=2 models can be coupled in a trivial way. We take one
sample of the percolation explorer and loop-erased that single sample in
two different ways to produce one sample, $\gamma^1$, of the gap=1 model
and one sample, $\gamma^2$, of the gap=2 model. One might expect that
if $d()$ is some distance function on curves, then 
for all $\epsilon>0$, $P(d(\gamma^1,\gamma^2) \ge \epsilon)$
converges to zero as the lattice spacing goes to zero. However, our
simulations provide strong evidence that this is not the case.

Since it appears that the scaling limit of the loop-erased percolation
explorer is a conformally invariant process that is not any SLE$_\kappa$,
an obvious question is what is it? A natural conjecture is that it is
some SLE($\kappa,\rho$) process. If the dimension of the curves is $4/3$
as the simulations indicate, then $\kappa$ would have to be $8/3$. But
we see no obvious conjecture for the other parameter(s). 

As for future work, one could consider loop-erasing other lattice models whose
scaling limit is SLE$_\kappa$ with $\kappa>4$, e.g., FK percolation,
and ask if the scaling limit is conformally invariant. One could also attempt
to define a loop-erasure process for SLE$_\kappa$ itself when $\kappa>4$.
The difficulty of course is that there is no natural linear chronological
order for the loops.

\bigskip
Acknowledgments:
This research was supported in part by NSF grant DMS-1500850.
An allocation of computer time from UA Research Computing 
at the University of Arizona is gratefully acknowledged.

\end{document}